\documentclass[onecolumn,
		nofootinbib
		]{revtex4}

\usepackage[utf8]{inputenc} 

\usepackage{graphicx}

\begin{document}

\title{The classical origin of modern mathematics}


\author{F. Gargiulo$^{1}$, A. Caen$^{2}$, R. Lambiotte$^{1}$, T. Carletti$^{1}$}
\affiliation{
1. Department of Mathematics and Namur Center for Complex Systems - naXys
\\ 2. DICE, Inria, Lyon, France}

\begin{abstract} 
The aim of this paper is to study the historical evolution of mathematical thinking and its spatial spreading. To do so, we have collected and integrated data from different online academic datasets. In its final stage, the database includes a large number ($N\sim 200K$) of advisor-student relationships, with affiliations and keywords on their research topic, over several centuries, from the 14th century until today. We focus on two different topics, the evolving importance of countries and of the research disciplines over time. Moreover we study the database at three levels, its global statistics, the mesoscale networks connecting countries and disciplines, and the genealogical level. \\

\end{abstract}
\maketitle
\section*{Introduction}
The statistical analysis of scientific databases, including those of the American Physical Society,  Scopus, the arXiv and ISI web of Knowledge, has become increasingly popular in the complex systems community in recent years. Important contributions include the development of appropriate scientometric measures to evaluate the scientific impact of scientists, journals and academic institutions \cite{eigen,jj,geo,radicchi,radicchi2} and to predict the future success of authors \cite{barabasiImpact,Acuna} and papers \cite{Shen2014}.
In parallel, the structure of collaboration has attracted much attention, and collaboration networks have become a central example for the study of complex networks, thanks to the high quality and availability of the datasets \cite{newman}. From a dynamical point of view, different papers \cite{barabasiCareer,resPath} studied the mobility of researchers during their academic career, showing that the statistical properties of their mobility patterns are mainly determined by simple features, such as geographical distance, university rankings and cultural similarity. \\

Limitations of the aforementioned datasets include their relatively narrow time window extending, at best, over 100 years and the difficulty to disambiguate author names, and thus to distinguish career paths across time. The original motivation of this paper was to address these issues by performing an extended study of The Mathematics Geneaology Project,  a very large, curated genealogical academic corpus \cite{mathgen}.  The dataset, whose basic statistics have been analysed elsewhere \cite{tesi,meso}, extends over several centuries and contains  pieces of information allowing us to retrieve the direct genealogical mentor-student links, but also university affiliations at different points of a career and  and the research domains. 
Data from the same website have already been used to assess the role of mentorship on scientific productivity \cite{amamral}  and   to study the  prestige of university departments \cite{porter}.

Our main goal is to analyse the history of modern mathematics, through the processes of birth, death, fusion and fission of research fields across time and space. In particular, we focus on the temporal evolution of the roles and importance of  countries and of  disciplines, on the structure of ``scientific families" and on the impact of genealogy on the development of scientific paradigms. As is often the case when performing a data-driven analysis of historical facts \cite{bbhist,bbhist2}, the data set is expected to be incomplete and to present  biases, mainly for the more ancient data. In the present case, the website collects the data in two ways: a participative method, based on the spontaneous registration of scholars (who can also register their students and their mentors), and a curated method, based on historical facts and performed by the creators of the web site.\\
The presence of biases calls for the use of appropriate  statistical measures, in preference based on ranking instead of absolute measures. In this work, we have also 
 introduced data-mining methods to correct and enrich the data structure. A first contribution of this work is thus methodological, with the design of a methodological setup that could be applied to other systems.
 We have then performed an analysis of the system at  three levels of granularity. First, a global one investigates the fully aggregated ``demography" (population in terms of countries and disciplines) of the database, with the aim
 to classify countries and disciplines according to their normalized activity behavior.
 Tracking the evolution of the rankings helps identify
transition points in the mathematical history, associated to emerging fields of research.
 Second, we have constructed directed weighted networks where nodes are scholars endowed with a set of attributes (thesis defense date, thesis defense location, thesis disciplines) and linked to other nodes using the genealogy associated to the mentor-student relation. This ``mesoscale" network allows us to investigate the relationships between the attributes and to identify a strong hierarchical structure in the scientific production in terms of countries as well as its evolution in the course of time.
 Finally, using an approach typic of kinship networks studies \cite{white,socNet}, we  focus on the statistical properties of the tree structure of the genealogy in terms of family structures. We  conclude by showing the presence of strong memory effects in the network morphogenesis.

\section*{ Dataset and associated networks}
The core of our dataset has been extracted from the website "mathematical genealogy project". It is one of the largest academic genealogy available on the web, consisting of approximatively 200K not--isolated scientists (186505) with information on their mentors and students. The data cover a period between the 14th century until nowadays. For a majority of mathematicians, we have detailed information about his/her PhD, including the title (for the 88\% of the scholars), the classification according to the 93 classes proposed by the American Mathematical Society\cite{ams} (for the 43\% of the scholars), the University delivering the degree as well as the year of its defense. However, because a large part of the database is spontaneously filled by the scientists, the data is imperfect and  attributes may be wrong or missing. A first step has thus consisted in comparing the database  with additional data  form Wikipedia  \cite{wiki} and, for more recent entries, with the Scopus profiles of scientists \cite{scopus}. After this preliminary phase, we have enriched  the  information available for the authors, by developing algorithms aimed at correcting the dates  based on the genealogical structure and the available statistics, assigning to each thesis a discipline,  based on the thesis title, and by disambiguating universities names. As previously stated, the algorithms, summarized  in the supplementary material, are  general  and could be applied in other contexts. After the enrichment, all the scholars of the database have a corrected date, the 88\% has an associated discipline and the 94\% and associated country. As a next step, we have exploited the enriched database  in order to study the geographical and temporal evolution of mathematical science. Different representations, described below, have been adopted.

\subsection*{The mesoscale networks}
\begin{figure}[h!]
\includegraphics[width=0.95\textwidth]{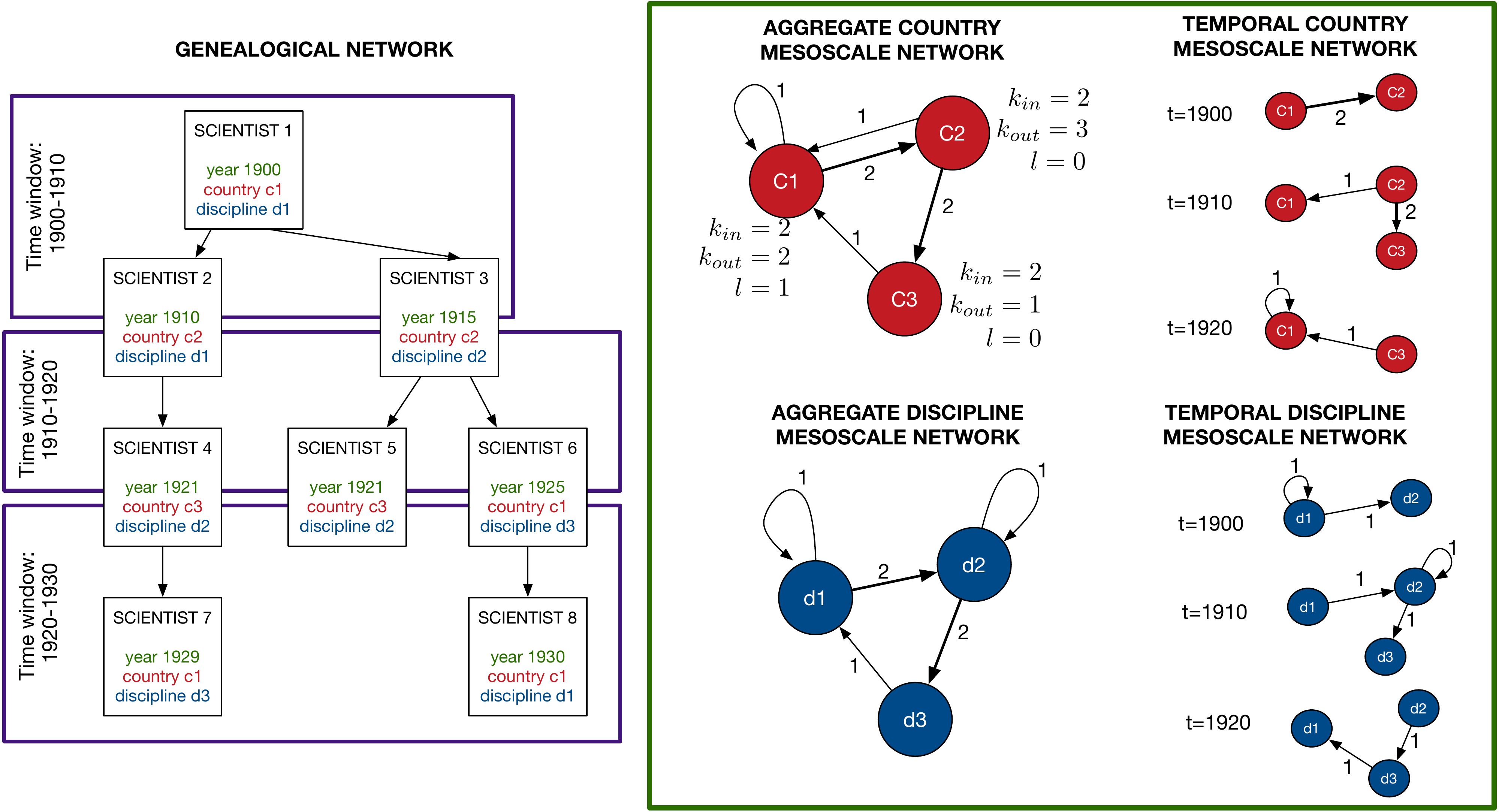}
\caption{\label{Fig0}
An example of the procedure to derive the mesoscale network from the genealogical data.}
\end{figure}

As a first step, we have considered a network where the nodes are  attributes of researchers (i.e.  universities, cities, countries and  discipline), and directed links, from $A$ to $B$, correspond to situations when a scientist with attribute $A$ was the PhD supervisor of a student with attribute $B$. In the case of universities, and under the assumption that a supervisor is in the same university as his/her PhD student, directed links are therefore a proxy for the mobility of a scholar, but also of the flow of knowledge between different places. In the case of disciplines, directed links correspond to a transfers of knowledge from one scientific discipline (the one of the mentor) to another one (the one of the student), from one generation to another, and how disciplines at a certain time may inherit, in terms of ideas and methods, from research fields at previous times. Let us note that the network allows for self-loops.
The procedure is illustrated, for the case of flows between countries, in Fig~\ref{Fig0}. 

\section*{Figures}

In the following, we will perform a  longitudinal study of the system, by considering the evolution of networks observed in different time windows. Note here that 
 the data are not uniformly distributed across time, with a strong bias towards recent times.

\subsection*{The genealogical tree and its partitions into families}
The genealogical graph is the most obvious representation of our dataset, consisting in an oriented acyclic graph \cite{DAG} linking a mentor to her/his students. This defines automatically the structure of hierarchical generations. Notice however that the structure of our data is not simply a tree due to the several cases where a student has two advisors.
A very common process in kinship is to cut the genealogical directed acyclic graphs into linear trees (alliances) where each individual has a single progenitor (the mother for representing the uterine links and the father for the agnatic ones). In this representation, the links between alliances represent the matrimonial structures between the different alliances in the society. In our context, when a scientist has more than one advisor, it is not clear which links should be cut to retrieve the original ancestors (no gender--like roles exist in this framework and moreover our dataset prevents us from identifying the principal supervisor from the secondary one, if any).  We thus propose a method to reproduce the optimal ancestry lines and to identify the important families in the genealogy. The method, fully described in the supplementary material, is based on the decomposition of the network into pure linear trees, and their statistical clustering based on probabilistic arguments; roughly speaking given two nodes $A$ and $B$ that can be linked in more than one way, thus implying the presence of non-trivial loops, we assign to every link in such paths the probability that $A$ and $B$ will be disconnected if the link is removed. We thus select links to be removed by maximising the probability that $A$ and $B$ are still linked. The resulting partition of the graph  into families identifies 84 families; remarkably, the $24$ mot populated families cover the $65\%$ of the scientific population in the database. Let us observe that alternative methods for family identification do exist, see for instance \cite{entropyTrees, TR}.

\section*{Results}
\subsection*{Global statistics}
\begin{figure}[h!]
\includegraphics[width=0.95\textwidth]{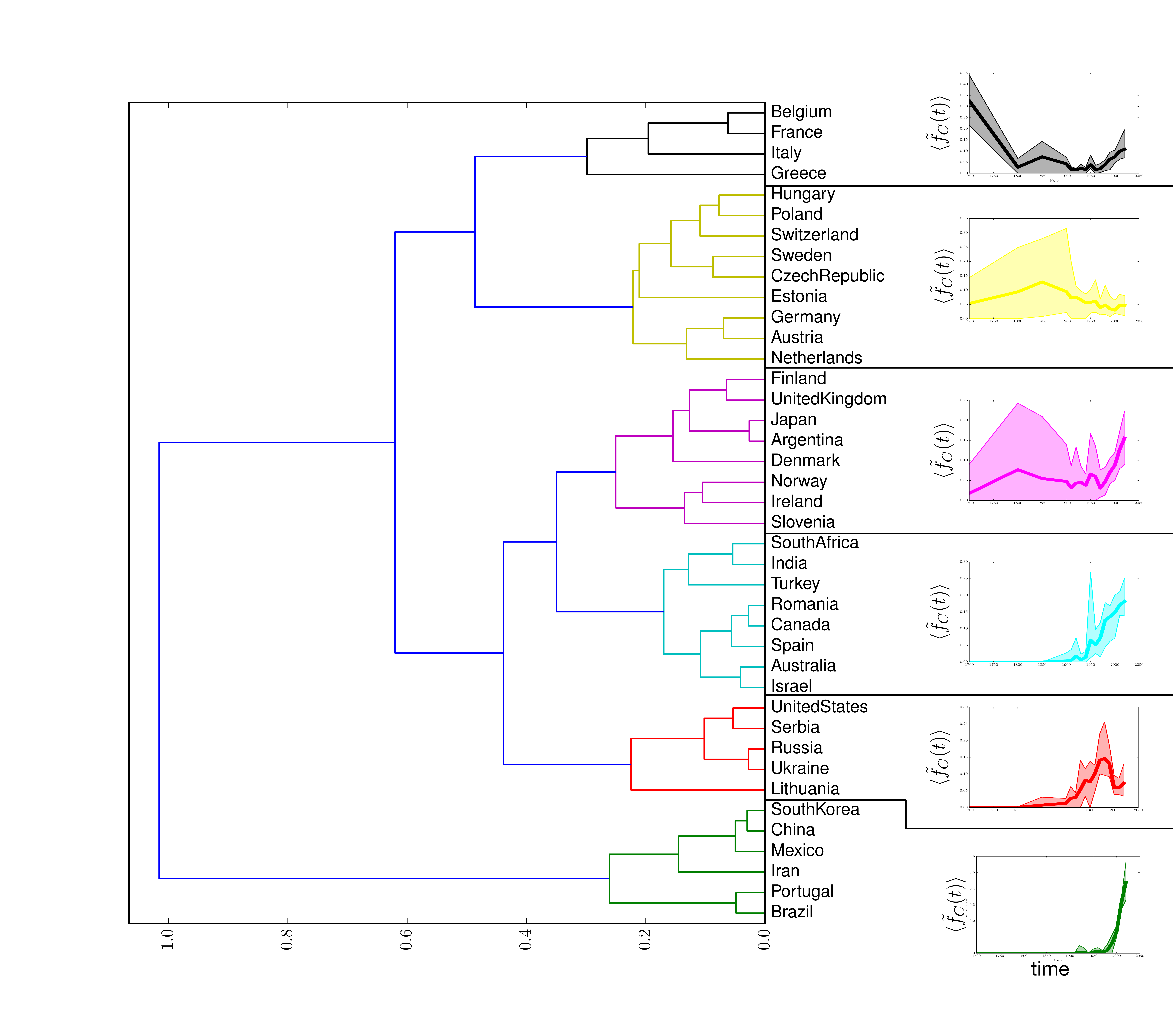}
\caption{\label{Fig1}
Clustering of the countries according to their time prevalence profile $\tilde{f}_I(t)$. The lines in the plots on the right describe the average prevalence profile for the cluster $\langle\tilde{f}_C(t)\rangle=\sum_{I\in C}\tilde{f}_I(t)/(\sum_{I\in C} 1)$. The shadowed area is included between the minimum value and the maximum value of the prevalence profile on the cluster ($\min_{I\in C}\tilde{f}_I(t)$,$\max_{I\in C}\tilde{f}_I(t)$ )}
\end{figure}

Let us define the {\em relative abundance profile} of each country in different periods $f_I(t)=N_I(t)/N(t)$, where $N(t)$ is the total number of scientists whose country is known in the database at time $t$ and $N_I(t)$ is the number of scientists in country $I$ at time $t$. Notice that, at this stage, the genealogical information is not used. From such profiles we can asses the evolution of the importance of the countries had in the history of mathematics due to the their different historical dynamics. In order to compare the profiles of different countries independently from their total production, we normalize each of these profiles by their \lq\lq volume\rq\rq ($\tilde{f}_I(t)=f_I(t)/\sum_tf_I(t)$) and then we classify them based on their Kolmogorov-Smirnov distance (see Fig.~\ref{Fig1} where the results are reported using a dendrogram, and the SI where we reported the profiles for the top 10 countries in the database). We observe different prototypical behaviors: countries with a central role in the ancient history whose centrality has decreased in the last centuries (for instance Italy, France and Greece), countries with a central role before the world wars (e.g. central Europe countries), countries emerging after the world wars (such as Japan and India), countries recently emerging (among which China and Brazil). Because of the normalization procedure we used, we obtain a cluster where USA is linked with ex-URSS countries and show a similar decreasing behavior in the latest decades (impossible to observed in a non--normalized context). \\

\begin{figure}[h!]
\includegraphics[width=0.95\textwidth]{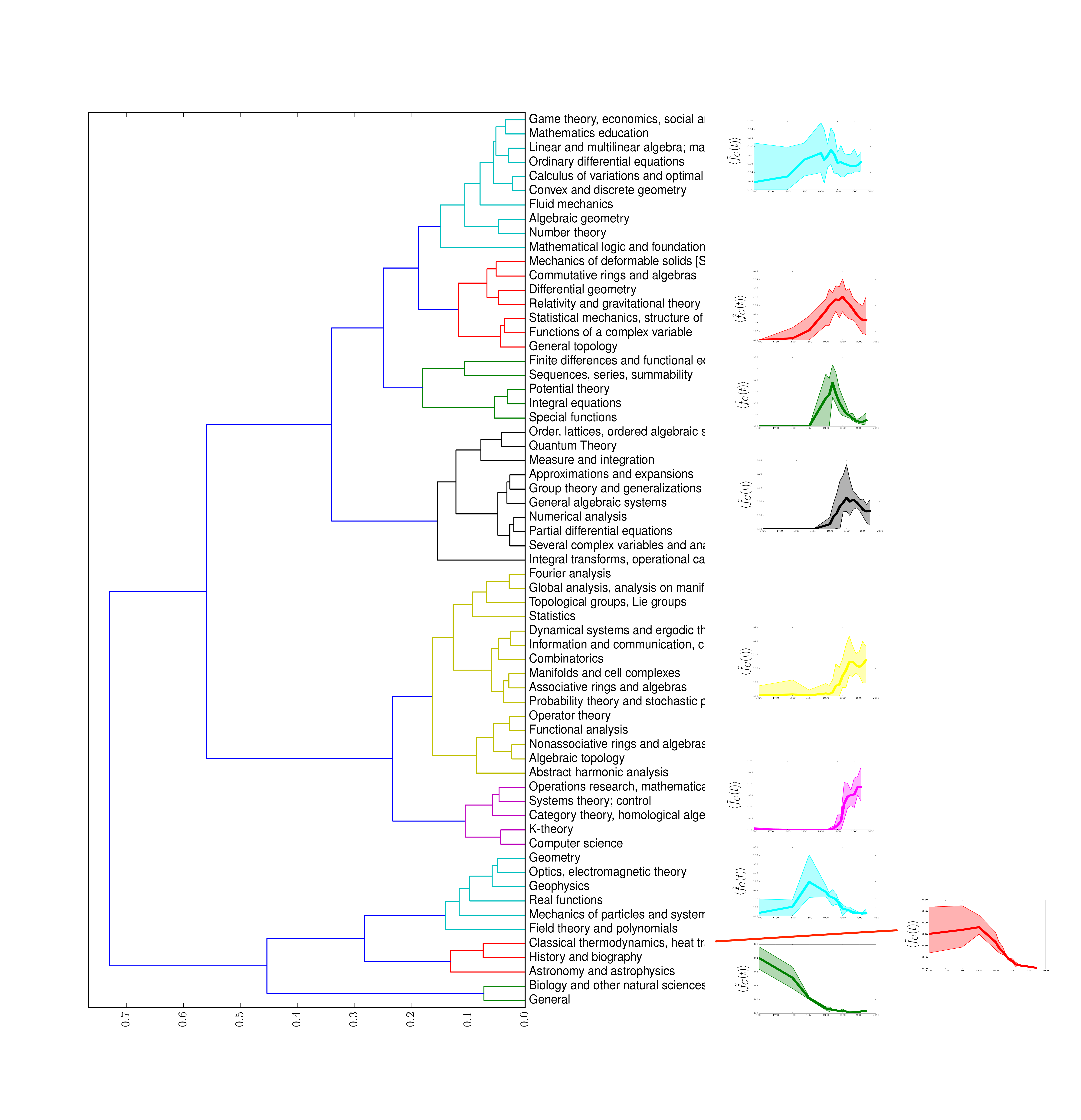}
\caption{\label{Fig2}
Clustering of the disciplines according to their time prevalence profile $\tilde{f}_I(t)$. The lines in the plots on the right describe the average prevalence profile for the cluster $\langle\tilde{f}_C(t)\rangle=\sum_{I\in C}\tilde{f}_I(t)/(\sum_{I\in C} 1)$. The shadowed area is included between the minimum value and the maximum value of the prevalence profile on the cluster ($\min_{I\in C}\tilde{f}_I(t)$,$\max_{I\in C}\tilde{f}_I(t)$ )}
\end{figure}

The same procedure is applied to disciplines and results reported in Fig.~\ref{Fig2} allow to identify three main blocks of disciplines: the disciplines that were more central during the industrial revolution (before 1900) are associated to physical applications (such as thermodynamics, mechanics and electromagnetism). The disciplines reaching their maximum of expansion around the 1950 are more abstract, even if several links exist to applied topics, such as telecommunication and quantum physics. Finally, the last decades have witnessed the emerging dominance of applied mathematics (e.g. statistics, probability) and computer science.  \\
To capture the rise and fall of countries or disciplines, we have compared the rankings of the top 10 countries and disciplines in different time periods. Standard indicators for rank comparison, such as the Kendall- Tau index, cannot be applied here since the elements in the  top-k lists are not conserved in time  \cite{kendalltau}. For this reason, we have used a distance measure based on a modified version of the Jaccard index allowing to compare ranked sets, $J(\textit{rank}_1,\textit{rank}_2)$ (more information is provided in the Supplementary material). As for the original Jaccard index the modified version is such that $J(\textit{rank}_1,\textit{rank}_2)=1$ when the rankings $\textit{rank}_1$ and $\textit{rank}_2$ are completely equivalent, and gives a value $0$ when these latter are not correlated at all. This information is then transformed in a distance by taking $d_J=1-J$. 
Increases in the distance measure, $d_J$, indicate major reshaping of the rankings.

\begin{figure}[h!]
\includegraphics[width=0.95\textwidth]{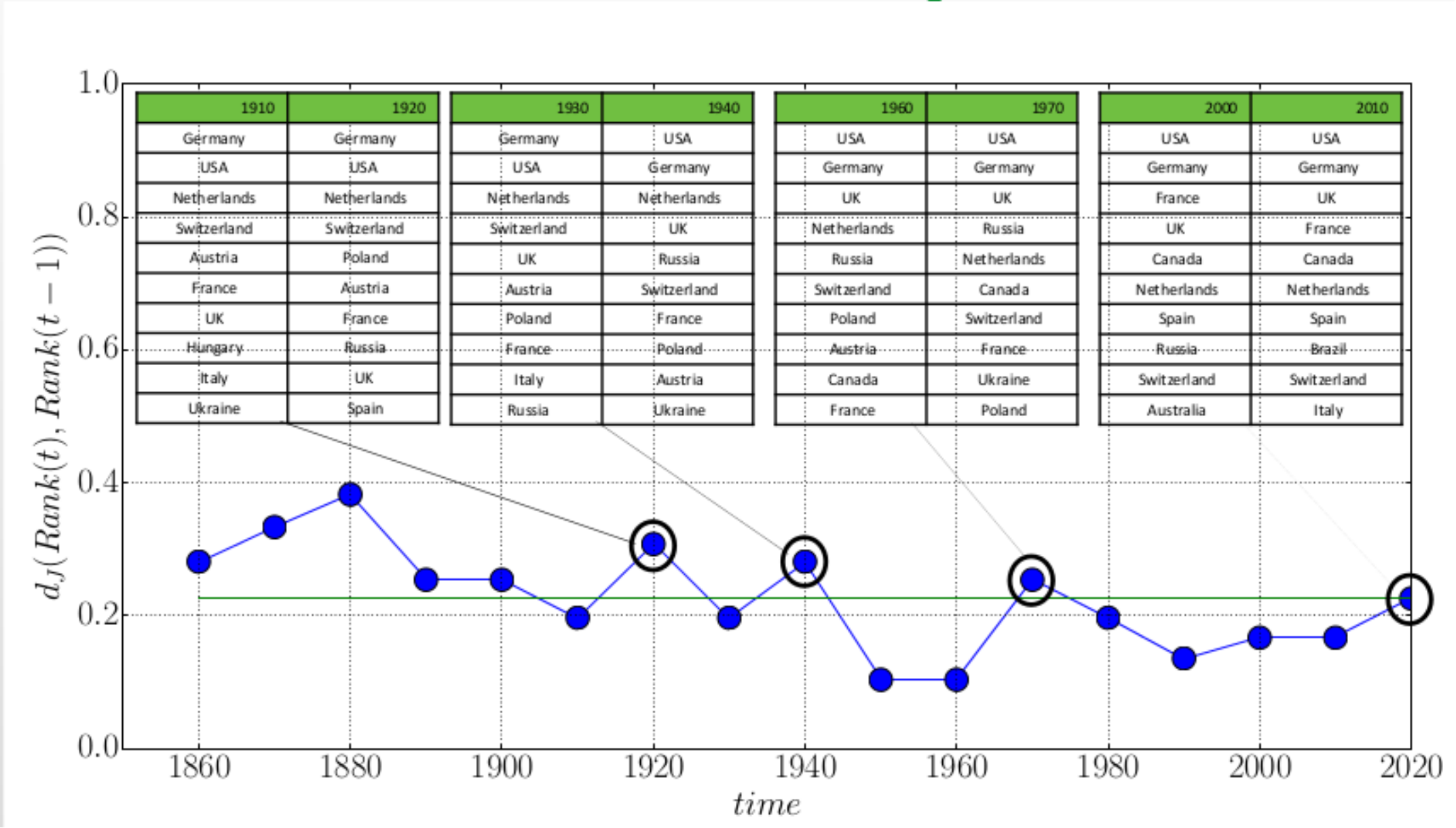}
\caption{\label{Fig3}
Modified Kendall-Tau index comparing the countries' rankings in different periods.}
\end{figure}

As we can observe in the upper plot of Fig.~\ref{Fig3}, corresponding to countries, we observe several transition points; for example a transition can be observed during the first World War, with the decreasing centrality of Austria and Hungary due to the end of the Austro-Hungarian emperor and the entering in the ranking of Russia.  Another transition is connected with the European political reshaping during the second World War, this is the period at which for the first time, USA surpasses Germany in the ranking. A third transition, around the 1960s, shows the increase of centrality of the Soviet Union (testified by the presence of several east-european countries in the ranking). Finally, more recently, we observe the decline of Russia and the emergence of new countries such as Brazil. \\

\begin{figure}[h!]
\includegraphics[width=0.95\textwidth]{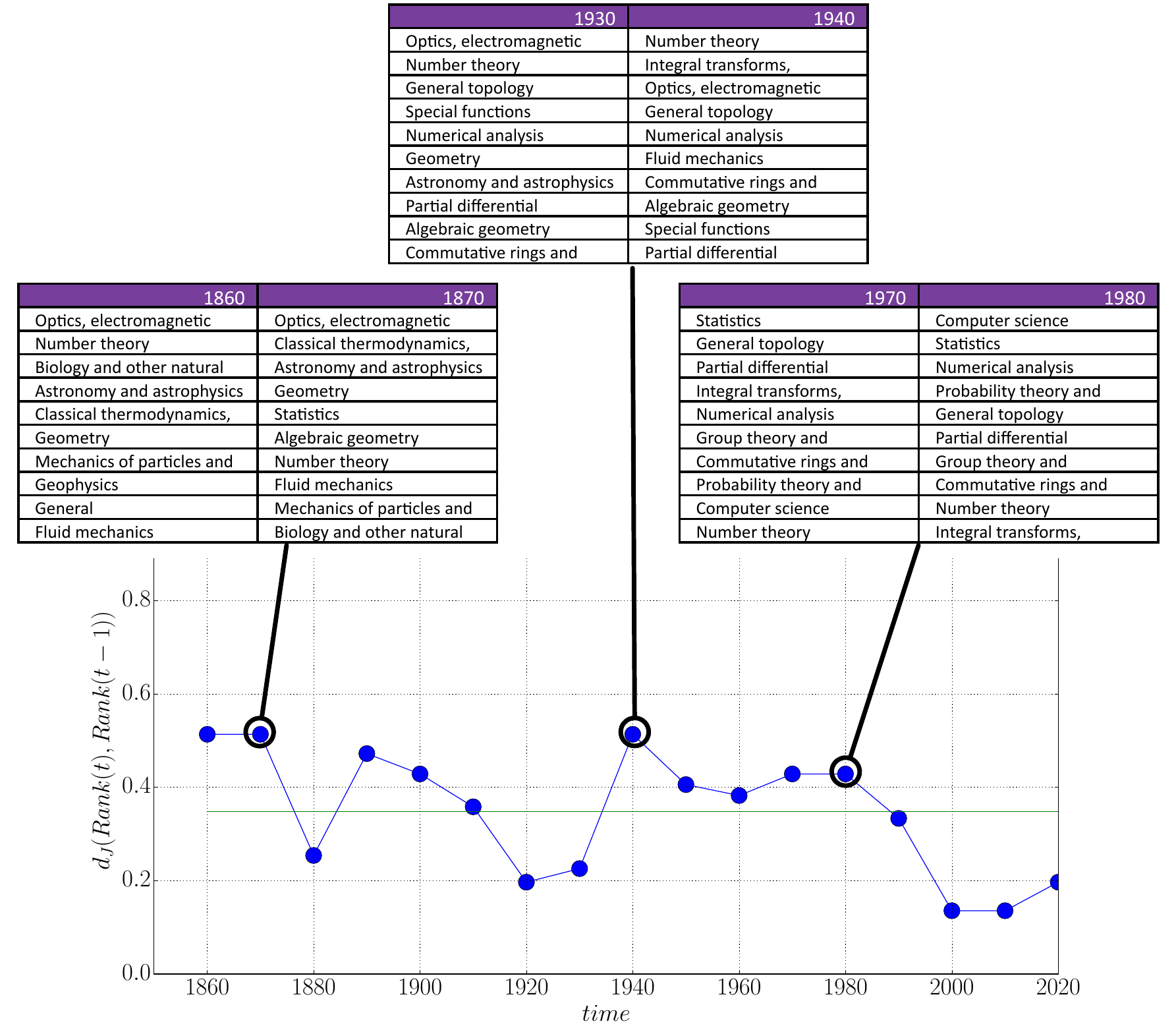}
\caption{\label{Fig4}
Modified Kendall-Tau index comparing the disciplines' rankings in different periods.}
\end{figure}

A similar analysis can be performed for disciplines. Results reported in Fig.~\ref{Fig4} show the presence of three tipping points. The first one is connected to industrial revolution and to the emergence of disciplines related to the physics of machines (such as thermodynamics and electromagnetism). The second one is connected to the emergence of fields linked to telecommunication and cryptography (e.g. number theory, spectral functions) during the second World War period. Finally the third one, in the 80's, concerns the emergence of computer science and statistics. 

\subsection*{Mesoscale networks}
\subsubsection*{Network of countries}
The countries network can be used to represent the {\em knowledge flows} from one country to another one, associated to the transition of a student in a country, becoming a professor in another country. The network presents few important hubs, that are the gravity centers of the scientific research (USA, Germany, Russia, UK). Each of these hubs tends to be surrounded by a community of countries. These communities can be associated to historical divisions, for instance a large block connected to USA scientific production, the Commonwealth nations, the ex-Sovietic block, the central European countries. The betweenness of countries allows to detect countries at the infercace between different communities, such as France connecting the central European countries with the USA-centered community or Poland connecting European research and the ex-Sovietic area.  \\ 

\begin{figure}[h!]
\includegraphics[width=0.95\textwidth]{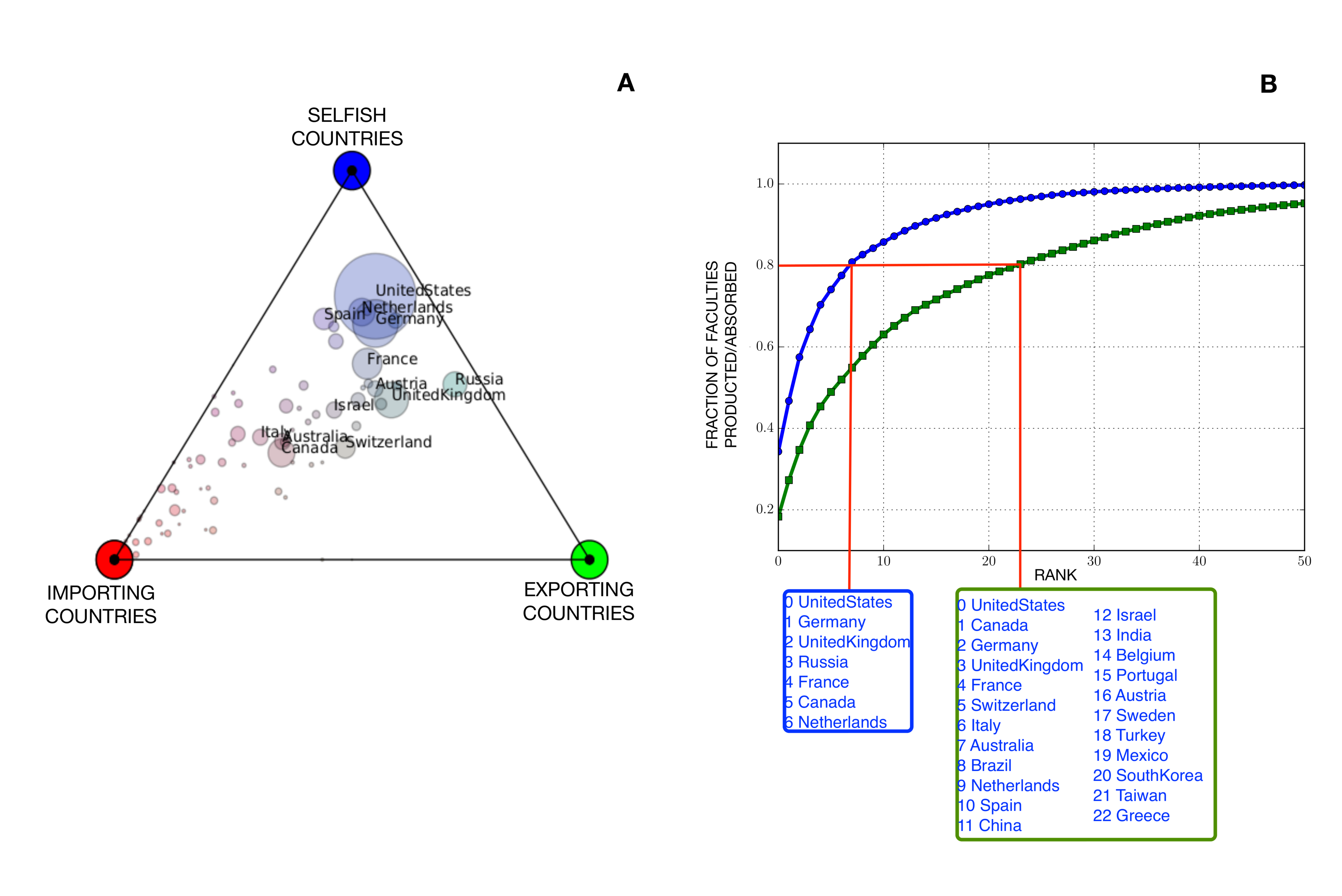}
\caption{\label{Fig5}Panel A: Relative position of the countries between selfish/exporting and absorbing behaviors. Panel B: Fraction of scientists produced and absorbed from the first countries in the rankings respectively of $k_{in}$ and $k_{out}$. In the boxes are displayed the countries producing and absorbing the 80\% of scientists. Both the panels concern the temporal aggregate of the network.}
\end{figure}

Another important index of the countries network is the weighted \emph{in(out)}--degree of the nodes and the number of self loops. The \emph{in}--degree of a country represents the number of scientists obtaining their PhD elsewhere and mentoring a PhD student in that country. Therefore a high \emph{in}--degree is associated to a country with a strong capacity of attracting scholars and absorbing knowledge from abroad. On the contrary, a high \emph{out}--degree represents a country producing scholars and exporting knowledge elsewhere. Each country can be therefore characterized by three normalized quantities, the fraction of scientist formed inside the country and remaining there, the fraction of scientist formed inside and leaving and the fraction formed abroad and absorbed by the country. In Fig.~\ref{Fig5}A we display the different positions of the most important countries with respect to these indexes, the closer a country is to one of the triangle vertices the larger is the associated index. The size of the dot  is a measure of the production of a country, i.e. estimated by its number of PhDs. The most productive countries tend to be the most \emph{selfish} ones, with a large fraction of self-loops. The most important \emph{exporters} are Russia and the UK. Countries with small scientific production show a tendency for \emph{importing} scholars. Observe that these indexes evolve in time  (see Fig. 3 of the SI). An important inversion point between $k_{in}(t)$ and the $k_{out}(t)$ is often observed around the second World War. Moreover, a key signature of emerging scientific countries seems to be the presence of $k_{in}(t)> k_{out}(t)$. \\

\begin{figure}[h!]
\includegraphics[width=0.95\textwidth]{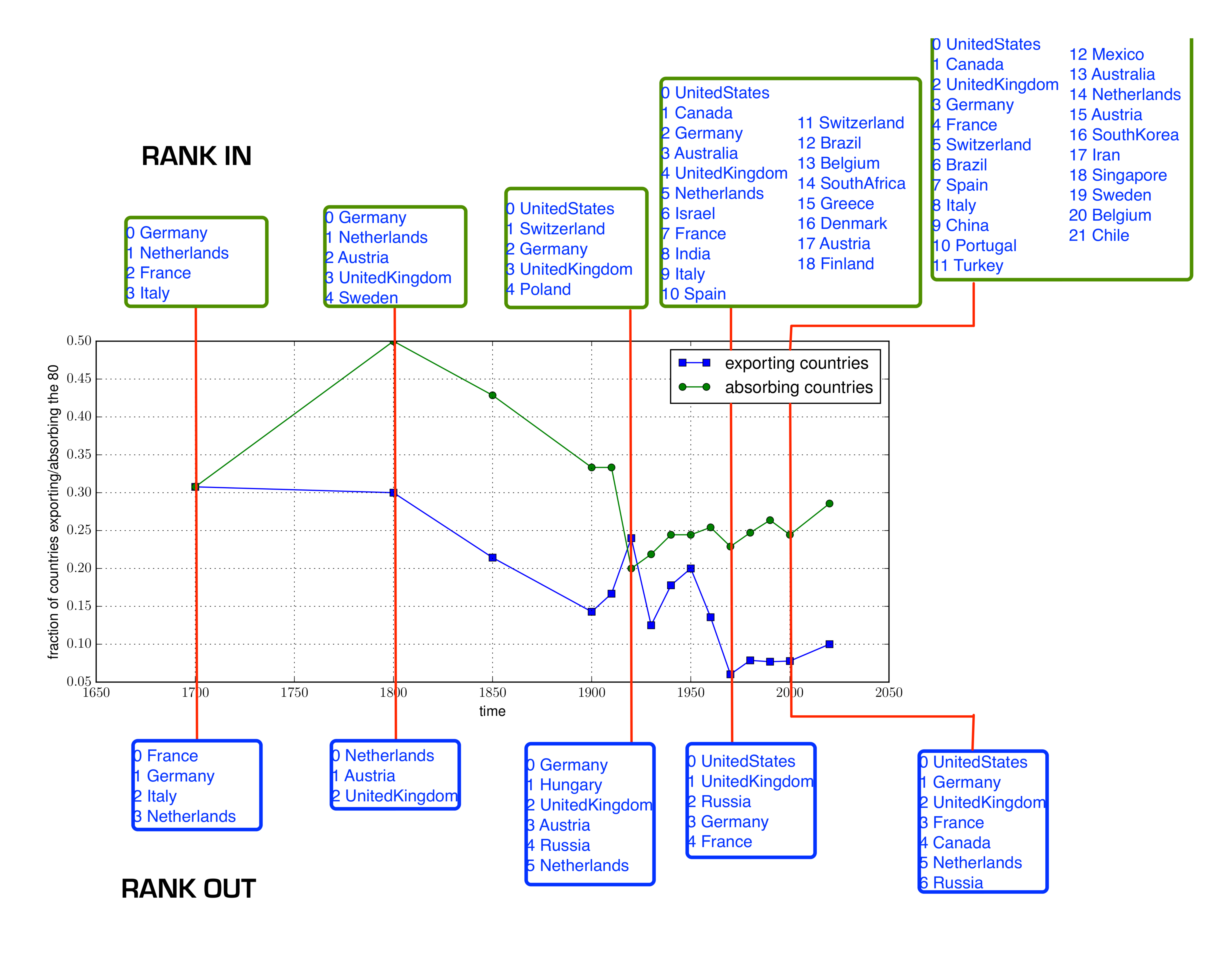}
\caption{\label{Fig6}Temporal network. Fraction of countries producing (and absorbing) the 80\% of the scientists in different historical periods.}
\end{figure}

To characterize the mobility of scientists across countries, we show in Fig.~\ref{Fig5}B the fraction of the total production and the total absorption of \emph{migrant} scientists for the first $r$ countries respectively in the \emph{in} and \emph{out}--degree ranking. The top 7 countries produce 80\% of the total international scholars. On the contrary, the curve concerning the  total absorption has a lower slope, depicting a larger worldwide spreading around the world. This shows a strong hierarchical structure in academic research where few countries ensure a large share of the worldwide diffusion of scientific knowledge. This scenario obviously evolved in time, as we observe in Fig.~\ref{Fig6}. Remark that scientific leaders changed at different points in time, but also that the scientific leadership group (countries producing the 80\% of the whole scientists production) is more restricted in recent times. It is interesting to notice that the minimal size of the scientific \emph{elite} has been reached in the sixties during the world bi--polarization resulting from the cold war. Since then, the size increased again with the emergence of globalization. 

More information about this network, in particular the properties of the aggregated transition networks concerning the whole historical period, can be found in the Supplementary Information.  

\subsubsection*{The transition network of disciplines}
The transition network of disciplines represents transfers of knowledge from one scientific discipline (the one of the mentor) to another one (the one of the student). The structure of this graph is quite homogeneous in terms of degree and four major communities can be identified: computer science, geometry, analysis and physics. Each community represents the disciplines exchanging more knowledge between them then with other research fields, and therefore can be interpreted as the scientific paradigms (according to Thomas Kuhn definition) at a certain period.\\

\begin{figure}[h!]
\includegraphics[width=0.95\textwidth]{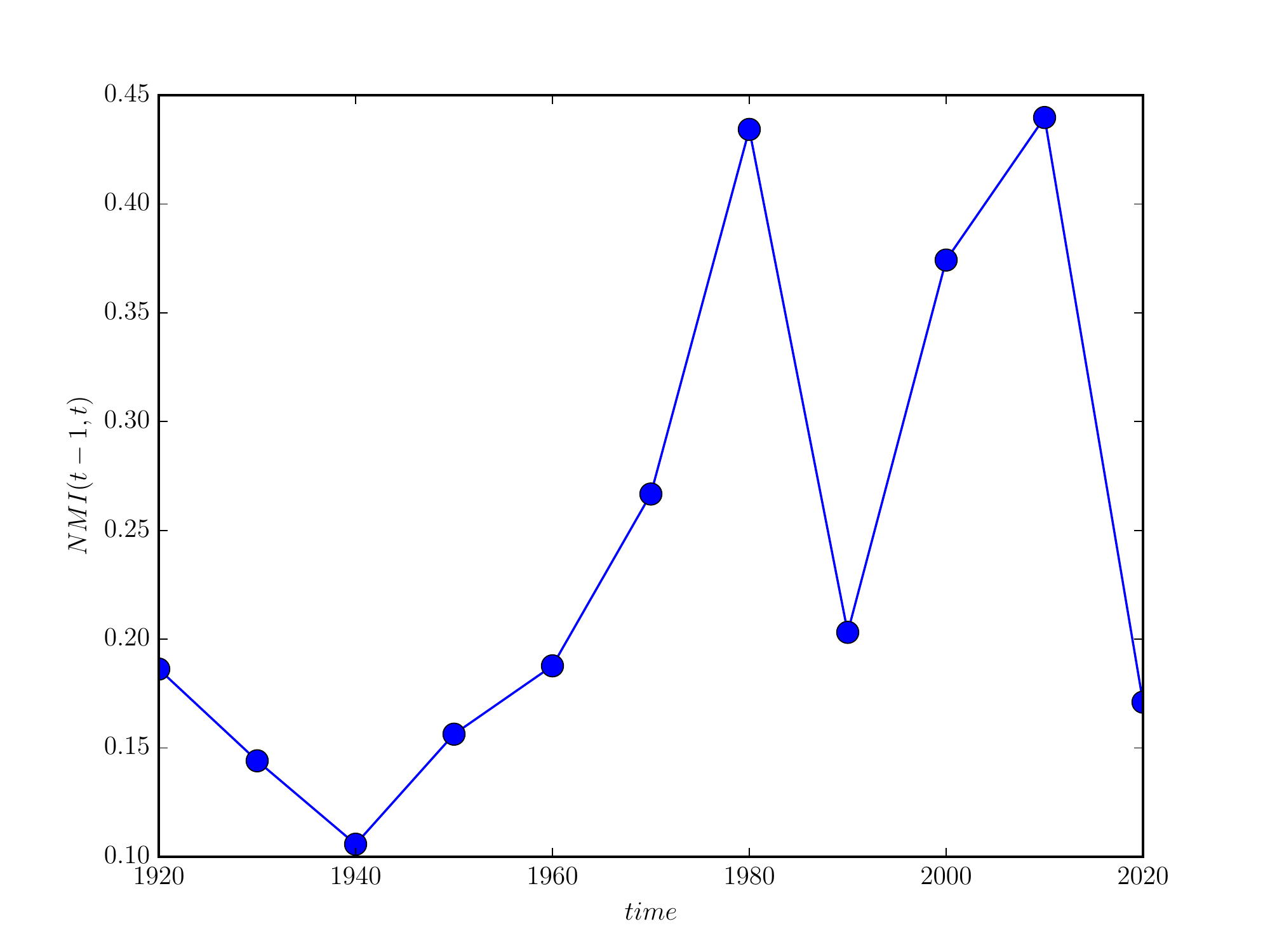}
\caption{\label{Fig7}Mutual information between the communities structures of the discipline network in different historical periods.}
\end{figure}

In Fig.~\ref{Fig7} we show the normalized mutual information (NMI) between the community structures obtained from different temporal slices of the network. The NMI index varies between one, when the two partitions are equal, and zero, when the two classifications are completely disjoined. A low value of the NMI indicates a ``revolution" in the sense of a strong reorganization of the knowledge structures, previously non interacting research fields start to exchange knowledge. The figure shows two important points where the NMI is low. The first transition, observed between 1930 and 1940, can be associated to the period when Statistics and Probability merged together, attracting then more applied disciplines like information theory, game theory and statistical mechanics, and leading to the emergence of the field of applied mathematics. The second transition is between 1970 and 1980, where computer science and statistics form one community, together with dynamical systems and applications in other fields of science. The latest transition is expected to be a spurious effect, due to a lack of data in recent years (the last time window starts in 2010 and therefore can contain data only for 5 years). Another potential approach, alternative to the measure of the NMI, to identify the structural changes in these structure could be the one proposed in \cite{change}.

\subsection*{The genealogical structure}
This last section is devoted to the study of the genealogy tree reconstructed from our data and of its relevance in the evolution of the history of mathematical science. \\

\begin{figure}[h!]
\includegraphics[width=0.95\textwidth]{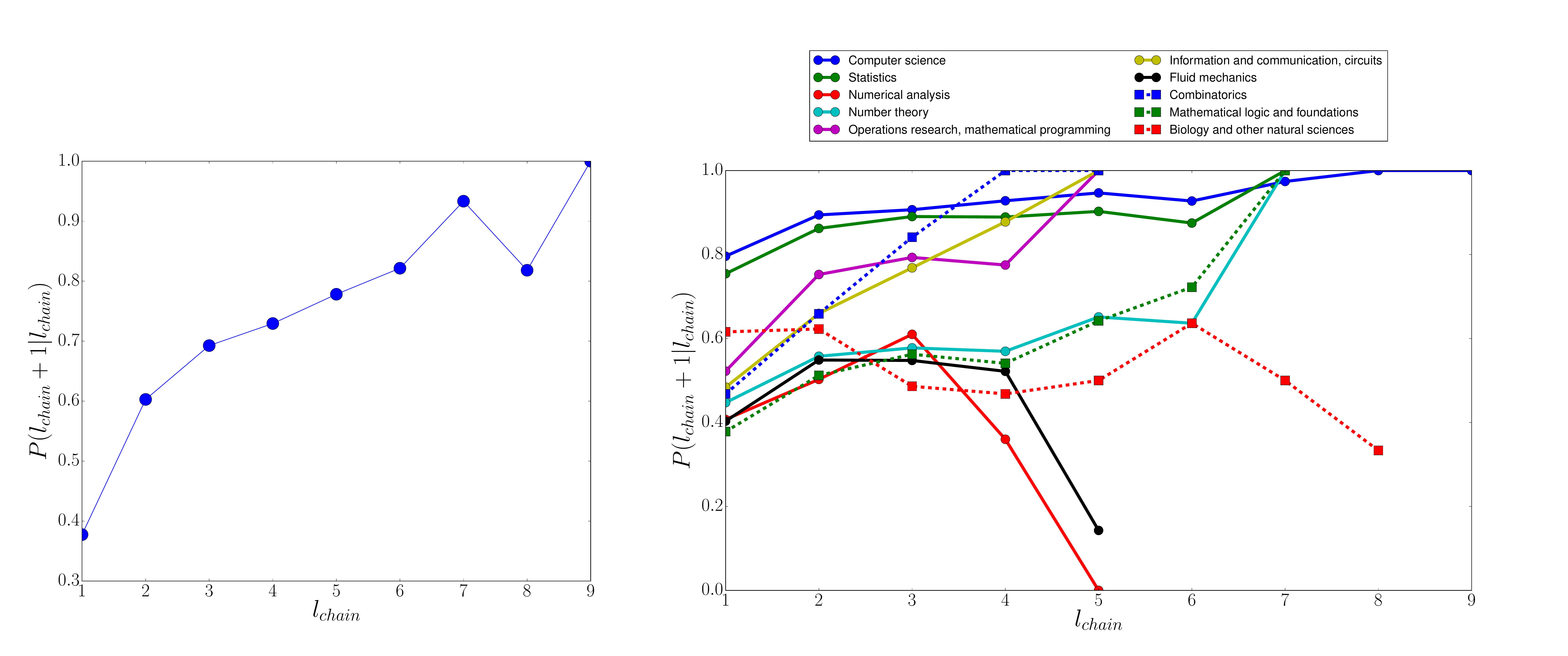}
\caption{\label{Fig8_v2}Conditional probability of having an iso--discipline chain of length $n+1$, having a chain of length $n$. Left panel: aggregated data all disciplines together; right panel: some selected disciplines.}
\end{figure}

The first result is the presence of a strong memory effects in the network morphogenesis, as students very often do research in the same discipline their mentor did. To quantify this idea we analyzed the genealogical chains where the ``filiation" link connects a mentor with a student maintaining the same discipline of the mentor. We call these objects iso-discipline chains. In Fig.~\ref{Fig8_v2} (left panel) we show results concerning the conditional probability of having a chain of length $n+1$ given a chain of length $n$, in other words the probability to have one more descendant working in the same research field of the whole chain, aggregating data over all the disciplines. The first point thus represents the probability for a student to have the same discipline of his/her mentor, one can clearly appreciate that as the chains get longer the probability to continue the same iso-discipline increases. This very marked memory effect in the network can be associated to the existence of ``schools" where a long tradition in a discipline exists such that new students are attracted and continue the tradition. Observe however (Fig.~\ref{Fig8_v2} right plot) that this phenomenon strongly depend on the discipline . \\

\begin{figure}[h!]
\includegraphics[width=0.95\textwidth]{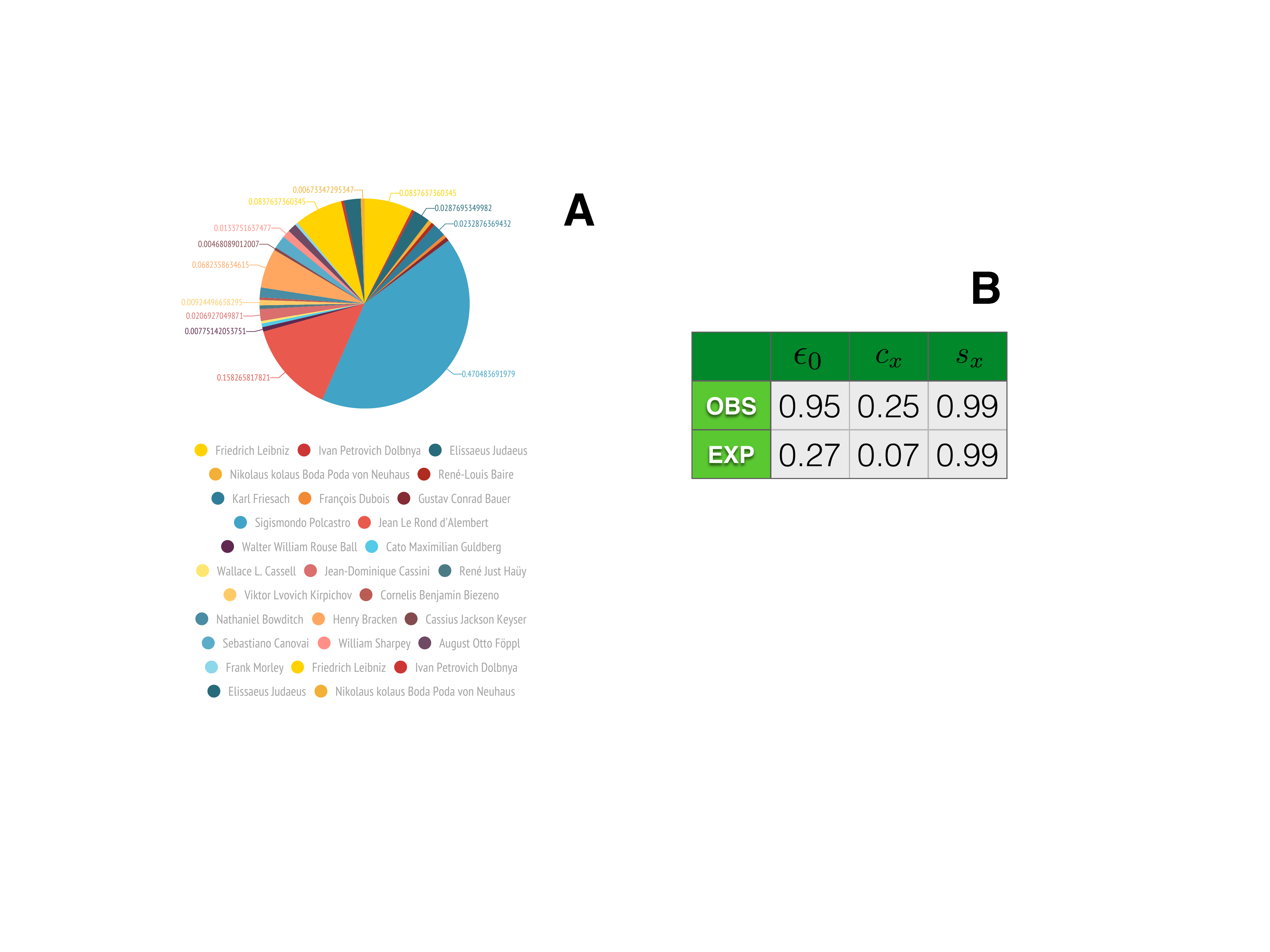}
\caption{\label{Fig9}Panel A: Relative size of the different families and family's initiator name. Panel B: Table with the values of the topological indicators for the real (observed) network in the first row and for the randomized model (expected) in the second row.}
\end{figure}

As previously explained, we have partitioned the network into disjoint families of scholars. Fig.\ref{Fig9}A shows that the 65\% of the scientists can be divided into 24 macroscopic families with size $S>500$. The largest family is the one originated in 1415 by the Italian medical doctor, Sigismondo Policastro. The second one, is the family originated by the Russian mathematician Ivan Petrovich Dolby, at the end of the 19th century. The large size of this family, born more recently than other families and geographically located mostly around Russia, is due to a high ``fecundity rate” in the Russian school of Mathematics.

The aggregated network between the families, reminiscent of kinship of the alliance networks defined in \cite{klaus,white}, can be described using some typical topological indicators~\cite{socNet}: 1) the endogamy index, $\epsilon_0$ describing the fraction of loops in the network (links between the same family); 2) The concentration index $c_x$ denoting the heterogeneity of the concentration of links between pairs of families ($c_x=1$ when all links are concentrated on a single pair and  $c_x=1/n^2$ when links are homogeneously distributed among the $n$ families); 3) the network symmetry index $s_x$ that varies from $0$ in case of total link unbalance, namely the outgoing flux and the ingoing one are very different each other, to $1$ in case of perfect symmetry of fluxes. To asses the relevance of such indicators computed for our genealogy network, we compared them with the expected values for a random multinomial reshuffling - null model - (see Fig.~\ref{Fig9}B), we can observe that, while the symmetry is a structural property, being unchanged by the reshuffling, the endogamy and the concentration are typical signatures of this network and moreover they are much higher than in traditional kinship networks ~\cite{socNet}. These results imply that the obtained scientific families are structurally very distant between them and that their relationships are very hierarchical (being these mediated by the largest families). \\

\begin{figure}[h!]
\includegraphics[width=0.95\textwidth]{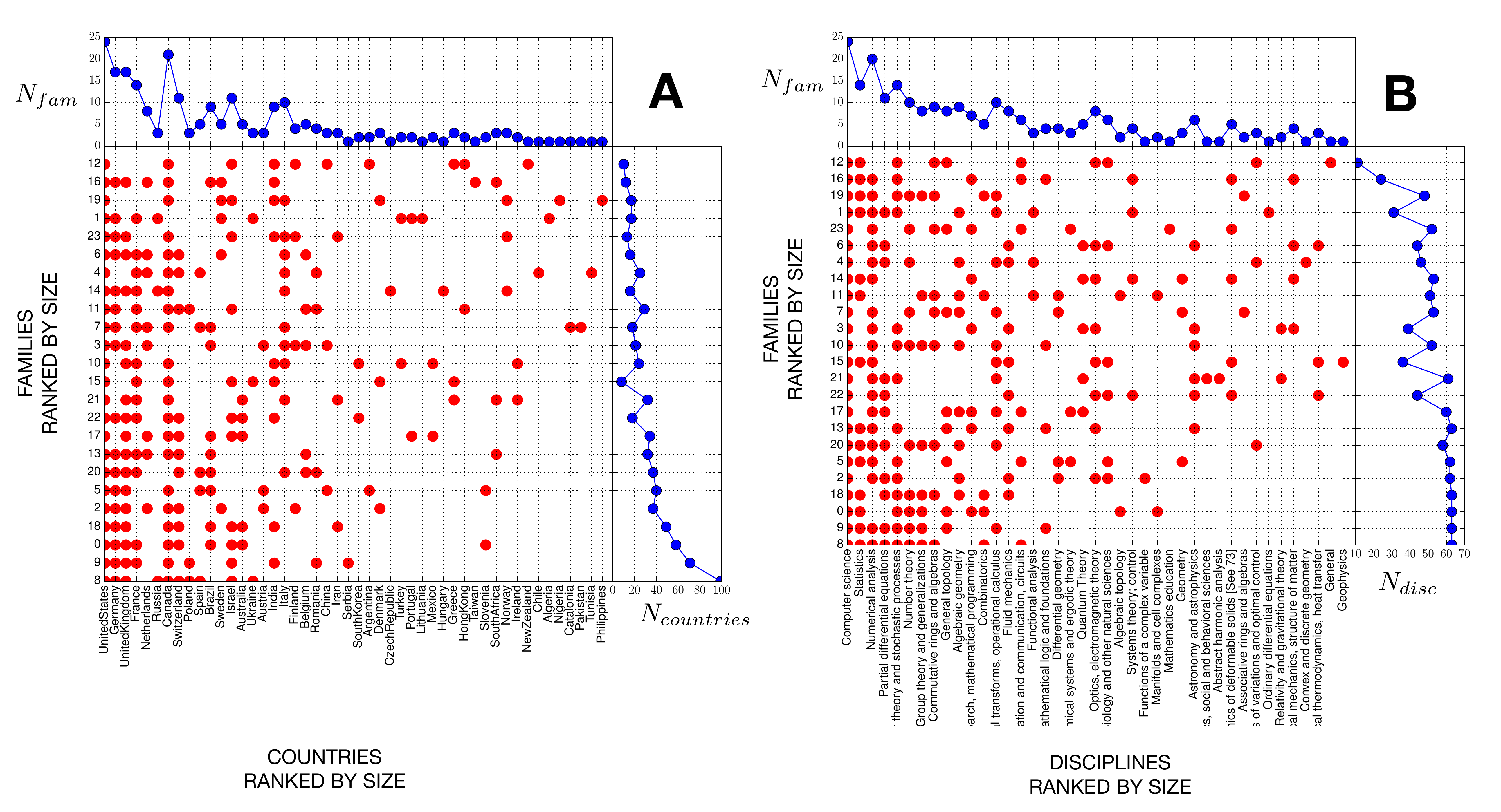}
\caption{\label{Fig10}Panel A: countries are set on the horizontal axis, ranked by the relative presence in the database, while in the  vertical axis, we report the families ranked by their size. A point at the intersection of the country-column and family-row indicates that the country is present in this family. The upper plot, is the column-marginal of the matrix represented in the central plot, representing the number of families where each country appears. The right plot is the row-marginal of the matrix represented in the central plot, representing  the number of countries present in each family. Panel  B: On the  horizontal axis we put the disciplines, ranked by the relative presence in the database, in the  vertical axis, the families ranked by the size. A point at the intersection of the discipline-column and family-row indicates that the discipline is present in that family. The upper plot, is the column-marginal of the matrix represented in the central plot, representing the number of families where each discipline appears. The right plot is the row-marginal of the matrix represented in the central plot, representing  the number of disciplines developed in each family.}
\end{figure}

Finally, we studied the distribution of families across countries and disciplines. As shown in Fig.~\ref{Fig10}A, with the exception of few cases, the most important countries (in term of production) are present in all the families, while the remaining countries are represented in a very low number of families (from 1 to 3). This feature implies a strong correlation between the genealogical structures and the geography. A similar behaviour can be observed for disciplines (Fig.~\ref{Fig10}B) even if, in this case, the curve describing the number of families with members working in a given discipline is smoother. We can therefore conclude that the genealogical families are strongly specialized in terms of geography and epistemic content.    

\section*{Conclusions}
In this paper, we have presented a data-driven study of the history of mathematical science, based on the Mathematical Genealogy Project. A first important aspect has been the cleaning and correction of the incomplete and sometimes inaccurate dataset. This operation was performed by means of machine-learning and by incorporating data from other sources, including Wikipedia.

We have then considered three different approaches to analyse the data: a demographic approach analyzing the time evolution of the prevalence of certain attributes (i.e. country or disciplines); a mesoscale network approach focusing on the connections between these attributes; a ``kinship" approach based on the clustering of genealogical trees. Our analysis reveals  important transition points in the history of mathematics and allows us to categorize countries according to their capacity to attract, export and self-maintain knowledge. Moreover, the community structures of the network of disciplines allows us to better describe the transformation of knowledge across time. Finally, we have  also identified important scientific families, associating them to their founder, and described their geographical and disciplinary distribution. 

Interesting lines of research for the future include the integration of additional datasets, based on different methodologies, to extend te scope of this work beyond the mathematical sciences. 

\bibliographystyle{bmc-mathphys} 


\appendix
\section{Methods}
In this document we describe the details of the statistical methods that we used in the paper.
\subsection{Rank comparison}
The Jaccard index is a statistical measure used to define the similarity between two sets $A$ and $B$:
\begin{equation}
J(A,B)=\frac{\# (A\bigcap B)}{\#(A\bigcup B)}\, ,
\end{equation}
where $\# A$ is the number of elements in the set $A$. By its very first definition it is thus a measure based simply on the content of the set, returning the fraction of elements that are common between two sets.
To measure the distance between ordered sets, on the other side it is often used the Kendall Tau distance, that is a measure of the number of permutations that happen between two rankings of the same set of elements.
In our case, the rankings in different time windows can contain (and it is usually the case) different sets of elements: a country, for example, can enter the ranking only at a certain time, being absent before. Therefore, in our case also the Kendall Tau index cannot be a suitable solution. \\
For this reason we introduced a new index allowing us to compare ordered sets with unequal entries, because it is based on the Jaccard index we decided to named it the "extended Jaccard", $\tilde{J}$. Given a set $A$ composed by $L$ elements and given an ordered classification of it, $r_A$, we define the unordered extended set $\tilde{A}$, of cardinality $\tilde{L}=\sum_{i=1}^Li=L(L+1)/2$ where each element of $A$ is repeated as many times as the complementary of its rank. For instance let $A=\{a,b,c,d\}$, thus $L=4$ and assume given the following ranking $r_A=\{1:a, 2:c, 3:d, 4:b\}$, that is $a$ is the first element, $c$ the second followed by $d$ and finally $b$, then the unordered extended set $\tilde{A}$ is given by:
\begin{equation}
\tilde{A}=\{a,a,a,a,c,c,c,d,d,b\}\, ,
\end{equation}
namely the set where each element $i\in A$ appears exactly $(L-r_i+1)$ times, being $r_i$ the rank of $i\in r_A$.
The modified Jaccard index is simply the Jaccard index applied to the unordered sets generated by this procedure:
\begin{equation}
\tilde{J}(r_A,r_B)=J(\tilde{A},\tilde{B})
\end{equation}
In such a way we can take into account the importance of a permutation with respect to the ranking of the exchanged elements and also consider new elements in the set because it can be directly applied also in the case where the entries of the rankings are not the same.

\subsection{Correcting the dates}
In the database we observed that $0.32\%$ of the couples mentor-student were characterized by an error in the reported dates, e.g. $year(mentor)>year(student)$. Moreover $6\%$ of the scientists are not associated at all with a date. To correct these errors and to fill the missing data we perform a stochastic iterative algorithm that modifies the dates in self consistent way using an heuristics method based on the distribution of the time distances between mentor-student and between students of the same mentor. The same procedure is also used to infer the missing data. 

\subsection{Associate a discipline to a scientist}
A fraction high as $88\%$ of the scientists presents the information about the title of the thesis. On the contrary only $43\%$ of the scholars has an associated disciplinary code (one of the $63$ classes of the AMS classification). We thus developed an algorithm able to infer from the title thesis, the most probable discipline. After a preliminary step where all thesis titles have been translated into english, we used a learning process based on the know thesis classifications to detect the main keyword in a title and the most probable discipline.

\subsection{Determine the families}
As already stated, in some case a student can have more than one mentor, in order to split the genealogical tree into disjoint families we developed an algorithm, based on the topology of the network, able to assign to each node with more than one "parent" its more probable familiar membership. The basic idea of the method is to generate random binary trees (namely random cutting one of the parents links for each scientist) and to look, for each realization of the network, if the detached parent and the son still fall in the same connected family. Repeatedly applying this procedure on all the possible random trees, we can extract for each node with multiple parental links, the probability for each link that its removal will not influence the positioning of the removed parent in the same family of the son. 

Since this procedure is computationally heavy, we performed a hierarchical grouping of the nodes in super-nodes structures: the super-nodes are the connected subgraphs of the network with a tree structure, in other words a super-node is obtained merging together all the nodes of the subtree under scrutiny. The super-nodes are then connected with weighed links obtained summing the number of links connecting the nodes part of each super-nodes. Because the total number of links is highly reduced, the calculation time of the most probable links is now strongly reduced. 

\section{The historical patterns}
\begin{figure}[h]
\centering
\includegraphics[width=1\textwidth]{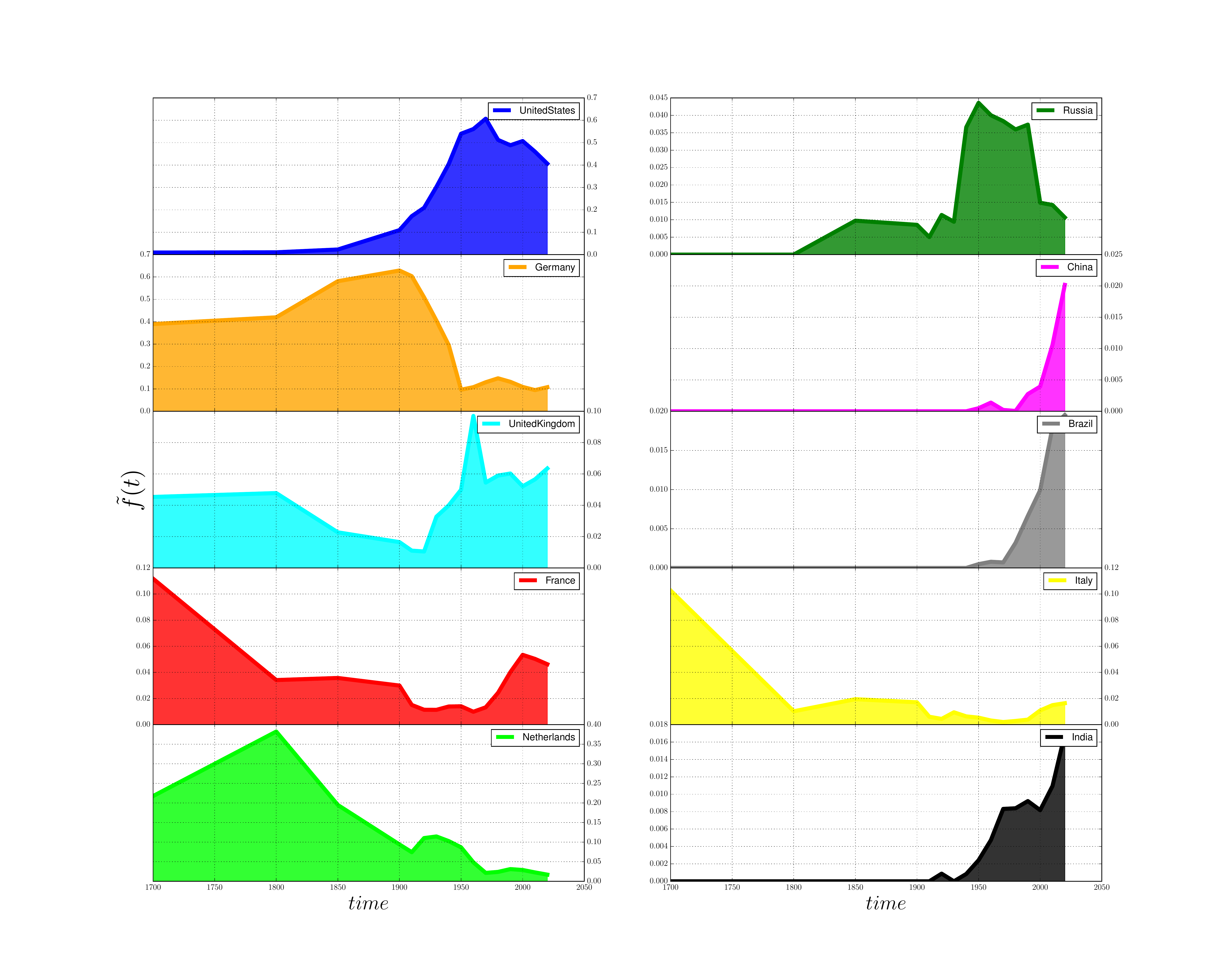}
\caption{\label{evolCount} Relative abundance profiles for the $10$ most significant countries in the database.}
\end{figure}

\begin{figure}[h]
\centering
\includegraphics[width=1\textwidth]{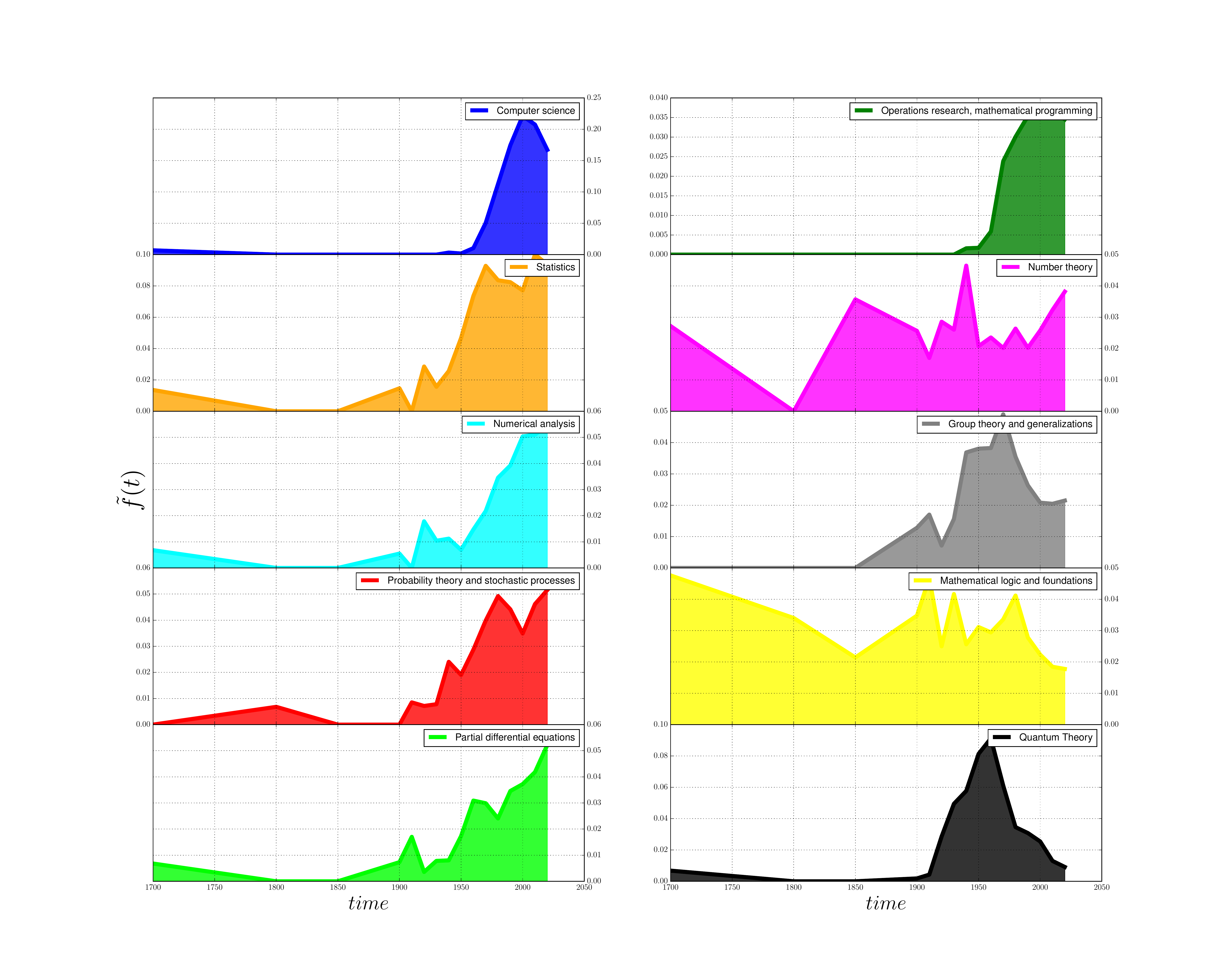}
\caption{\label{evolDisc} Relative abundance profiles for the $10$ most representative disciplines in the database.}
\end{figure}

In this section we show some prototypical "participation" patterns for different countries and disciplines, namely the profiles of the relative presence in each time window, for the $10$ most significant countries (see Fig.~\ref{evolCount}) and disciplines (see Fig.~\ref{evolDisc}). \\

From Fig.~\ref{evolCount} we can clearly notice the loss of centrality of countries like Italy and Netherlands, and partially France. The rise of US and Russia coincide with the decline of Germany after the Second World War. Finally we can evince the fast growing trends for emerging countries like China, India and Brazil. \\
In the case of disciplines, one can appreciate from Fig.~\ref{evolDisc} the growth and fall of quantum theory and the parallel growth of group theory that is strongly inherent the physics of the 60's. Other disciplines like number theory and logic are quite uniformly distributed in the time.

\section{The mesoscale network structures}
\begin{figure}[h]
\centering
\includegraphics[width=1\textwidth]{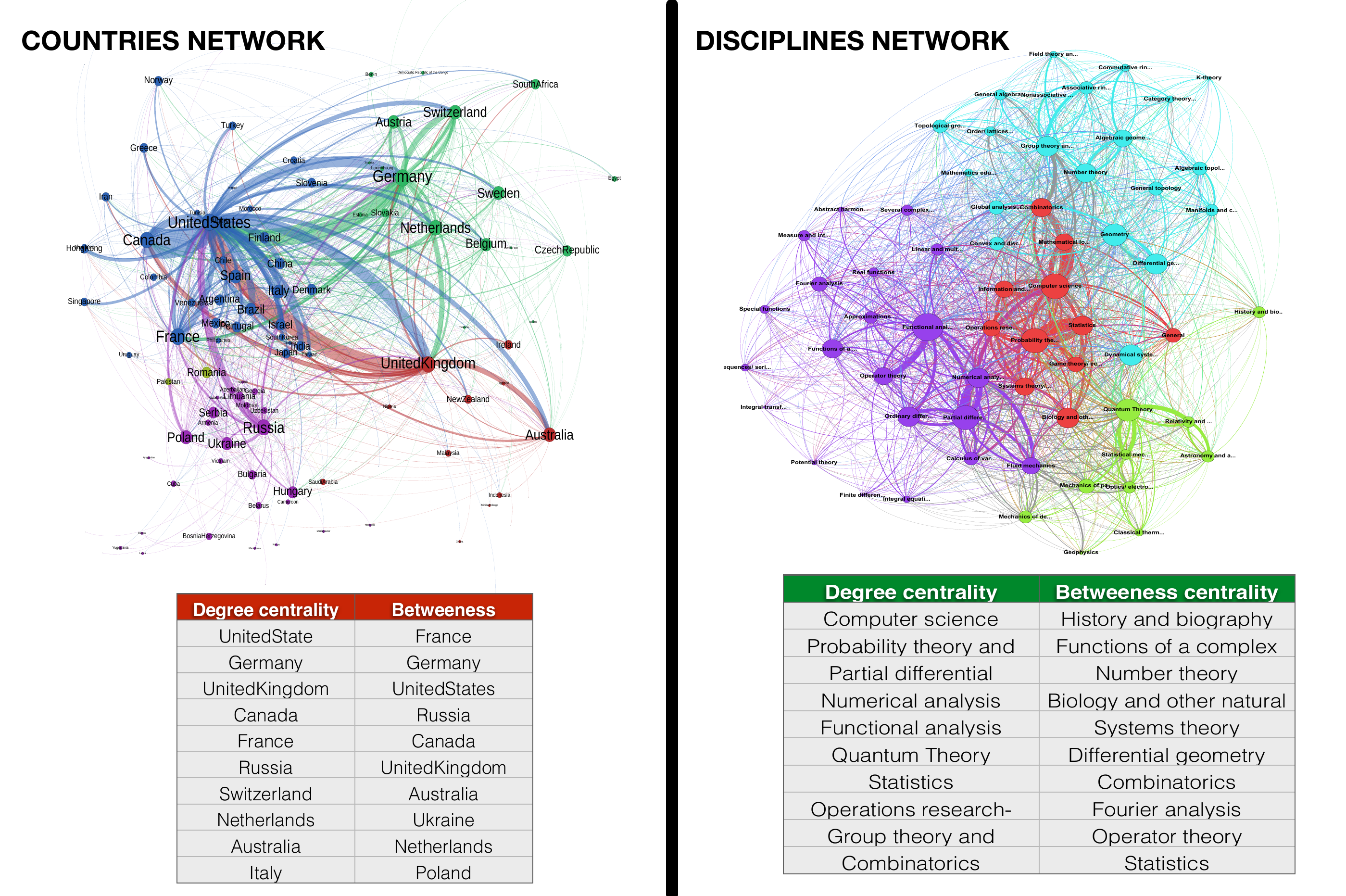}
\caption{\label{netSt} Static network representation for countries and disciplines. The node size is proportional to the total degree. Some centrality measures are shown in the tables. The colors of the nodes represent the community structures of the network determined using the Louvain algorithm.}
\end{figure}

In Fig.~\ref{netSt} we report the structure of the time-aggregated static mesoscale networks. As one can evince from the graph, the strengths of the country network are strongly heterogeneous, varying from $6500$ (USA) to $1$ for several peripheral counties appearing just once in the database. On the other hand the heterogeneity is less marked in the case of disciplines. As we can observe, both for countries and disciplines, the degree centrality top list is not equivalent to the betweenness. For the countries, the comparison between the top-10 countries in terms of degree and in term of betweenness put in evidence the very strategical role of France and est-European countries (Russia, Ukraine, Poland) in connecting information flows between the network. 

The difference between the two rankings is still more astonishing for the disciplines where just two disciplines (combinatorics and statistics) appear in the two top-10 lists. This scenario suggests the existence of very well structured epistemic communities around the top connected disciplines, where the connections between the communities is performed by peripherals nodes of each community. 

\begin{figure}[h]
\centering
\includegraphics[width=1\textwidth]{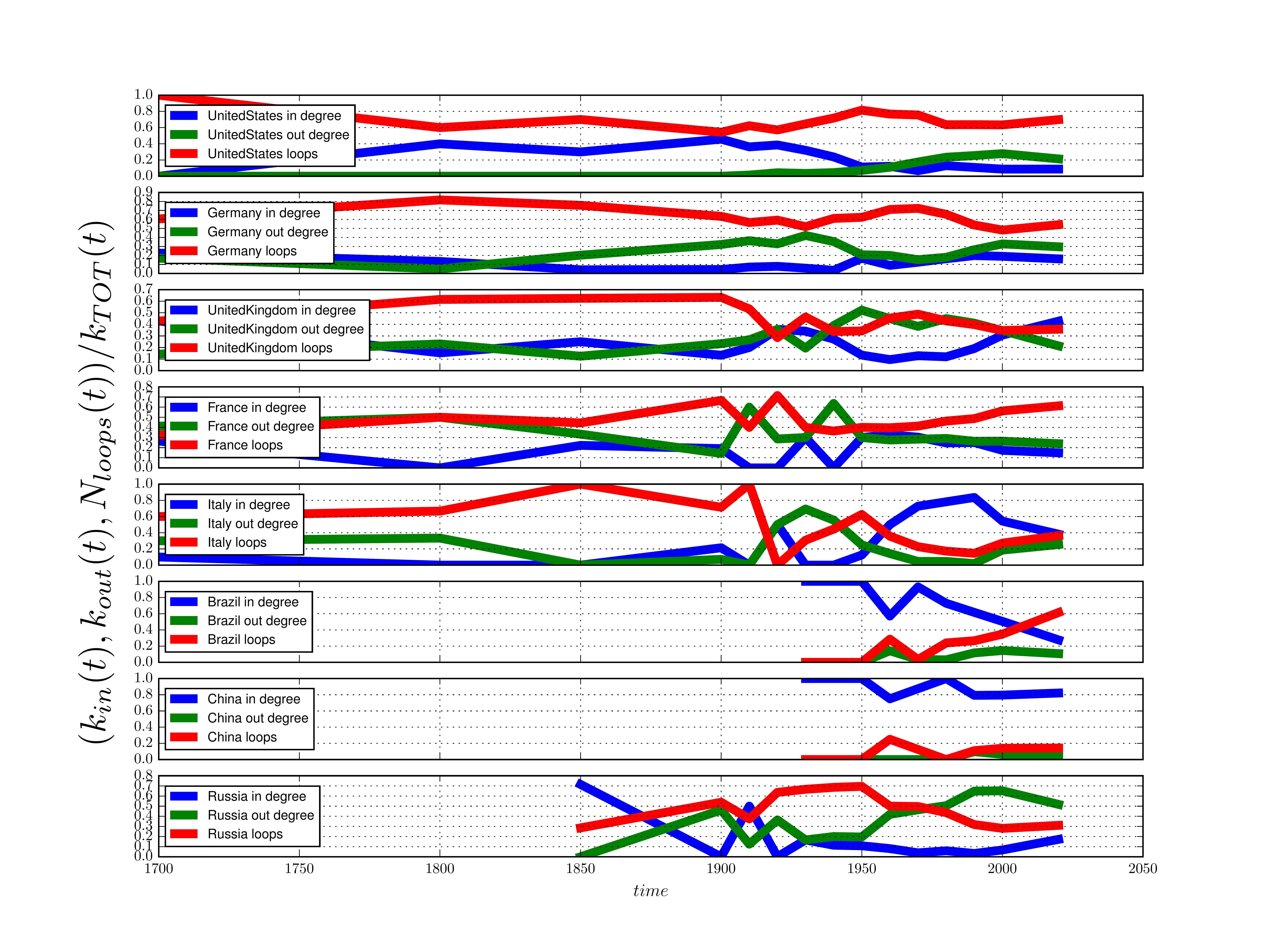}
\caption{\label{degreeEvol} Fraction of In-links,out-links and self-links for some representative countries in the dataset.}
\end{figure}

In Fig.~\ref{degreeEvol} we represent the time evolution of the fraction of self-loops, in- and out-links, with respect to the total number of links for a given country. As we can observe the emerging countries are characterized by a significantly high in-degree and a very small out-degree. The case of Italy is quite interesting showing a strong exodus during the second world war, followed by an inversion around the 50s due to the scientific politics of favoring the coming back of the emigrated scientists. The same maximum of the out-degree due to the second world war can be observed in Germany, while the opposite sign can be noticed in USA and UK. Notice also the following inversion of in and out degree prevalence for US after the 50s. This inversion is even more pronounced for Russia. 

\end{document}